\documentclass[12pt]{article} \usepackage{latexsym}
\usepackage{amsfonts} \usepackage{amsmath} \usepackage{amssymb}

\makeatletter
\@addtoreset{equation}{section}

\ifx\reset@font\undefined
           \let\reset@font=\relax
\fi
\def\section{\@startsection {section}{1}{\z@}{-3.5ex plus-1ex minus
             -.2ex}{2.3ex plus.2ex}{\reset@font\large\bf}}
\def\subsection{\@startsection{subsection}{2}{\z@}{-3.25ex plus-1ex
             minus-.2ex}{-.1em}{\reset@font\large\bf}}
\def\subsubsection{\@startsection{subsubsection}{3}{\z@}{-3.25ex plus
          -1ex minus-.2ex}{-.1em}{\reset@font\normalsize\bf}}
\makeatother

\title{Saturating Constructions for Normed Spaces}
\author{Stanislaw J. Szarek
\thanks{Supported in part by a grant from the
National Science Foundation (U.S.A.).}
                 \and  Nicole Tomczak-Jaegermann
\thanks{This author holds the Canada
    Research Chair in Geometric  Analysis.}
}

\newcommand\address{\noindent\leavevmode%
Equipe d'Analyse Fonctionnelle, BP 186\\
Universit\'e Pierre et Marie Curie\\
75252 Paris, France\\
{\it and}\\
Department of Mathematics\\
Case Western Reserve University\\
Cleveland, Ohio 44106-7058, U.S.A.\\
{\small\tt%
szarek@ccr.jussieu.fr}\\[.5cm]
\noindent
Department of Mathematical and Statistical Sciences,\\
University of Alberta,\\
Edmonton, Alberta, Canada T6G 2G1,\\
{\small\tt%
           nicole@ellpspace.math.ualberta.ca} }

\date{}

\newtheorem{fact}{Fact}[section]
\newtheorem{thm}[fact]{Theorem}
\newtheorem{prop}[fact]{Proposition}
\newtheorem{lemma}[fact]{Lemma}

\newbox\nrmbox
\setbox\nrmbox=\hbox{$\Vert$}
\def\nrmrule{\vrule height\ht\nrmbox depth1.2\dp\nrmbox}
\setbox\nrmbox=%
           \hbox{\kern0.15em\nrmrule\kern0.15em\nrmrule\kern0.15em%
                     \nrmrule\kern0.15em}
\newcommand{\Snorm}[1]%
           {\copy\nrmbox#1\copy\nrmbox\kern-0.03em%
                 \lower.4ex\hbox{}}

\newcommand{\rem}{\noindent{\bf Remark{\ \ }}}

\newcommand{\proof}{\noindent{\em Proof{\ \ }}}

\newcommand{\qed}{\bigskip\hfill\(\Box\)}

\newcommand{\R}{\mathbb{R}}

\newcommand{\N}{\mathbb{N}}

\newcommand{\pp}{\mathbb{P}}

\newcommand{\bt}{\beta}

\newcommand{\ep}{\varepsilon}

\newcommand{\ra}{{\rightarrow}}

\newcommand{\spn}{\mathop{\rm span\,}}

\newcommand{\conv}{\mathop{\rm conv\,}}

\newcommand{\esp}{{\tilde{E}}}
\newcommand{\kk}{{\tilde{K}}}
\newcommand{\dd}{{\tilde{D}}}

\newcommand{\wrt}{with respect to }

\begin{document}

\maketitle

{\abstract{\small
We prove several results of the following type:
given finite dimensional normed space $V$ there exists another
space $X$ with $\log \dim X = O(\log \dim V) $ and such that
every subspace (or quotient) of $X$, whose dimension is not ``too
small," contains a further  subspace isometric to $V$.  This sheds new
light on the structure of such large subspaces or quotients
(resp., large sections or projections of convex bodies) and allows to
solve several  problems stated in the 1980s by V. Milman.}}

\section{Introduction}
\label{introduction}

Much of geometric functional analysis revolves around the study of
the family of subspaces (or, dually, of quotients) of a given Banach
space.
On the one hand, one wants to detect some  possible regularities
in the structure of subspaces  which  might have not existed
in the whole space; on the other hand,
one tries to determine to what degree the structure of
the space can be recovered from the knowledge of its subspaces.
In the finite dimensional case there is a compelling geometric
interpretation: a normed space is determined by its unit ball,
a centrally symmetric convex body, subspaces correspond to sections
of that body,  and quotients to projections.
A seminal result in this
direction  was the 1961 Dvoretzky theorem: every convex body of (large)
dimension $n$ admits central sections which are approximately
ellipsoidal and whose dimension $k$ is at least of order $\log{n}$.
(A proof giving the logarithmic order  for  $k$ was 
given  by Milman in 1971 (\cite{M-dv}, or see \cite{MS1}.)
In other words, the most regular $k$--dimensional structure, i.e., the
Euclidean space, is present in every normed space of sufficiently high
dimension. Another interpretation of the same phenomenon is as
follows. From the information-theoretic point of view, the Euclidean
spaces have minimal complexity.  Thus by passing to a subspace or a
quotient -- of dimension much smaller than that of the original space,
but still a priori arbitrarily large -- one may, in a sense, lose all
the information about a space.

It has been known that the $\log{n}$ estimate in the Dvoretzky
theorem, while optimal in general, can be substantially improved for
some special classes of spaces/convex bodies.  Still, it was a major
surprise when V. Milman (\cite{M1}) discovered in the mid 1980's that
{\em every } $n$-dimensional normed space admits a {\em subspace of a
  quotient } which is ``nearly" Euclidean and whose dimension is $\ge
\theta n$, where $\theta \in (0,1)$ is arbitrary (of course, the exact
meaning of ``nearly" depends on $\theta$).  Moreover,
as the first step of his approach he proved that every $n$-dimensional
normed space admits
a ``proportional dimensional'' quotient of bounded volume ratio,
a volumetric  characteristic of a body
closely related to cotype   properties (we refer to \cite{MS1}, \cite{T}
and particularly \cite{P} for definitions of these and other basic
notions and results that are relevant here).
This showed that one can get a very essential regularity in a
global invariant of a space by passing to a
quotient or a subspace of dimension, say, approximately $n/2$ (while
of course also losing a substantial part of the local information).
It was thus natural to expect that similar statements may
be true for other related characteristics. This line of thinking
was exemplified in a series of problems
posed by Milman in his 1986 ICM Berkeley lecture \cite{M2}.
A positive answer to any of them would have many important
consequences in geometric functional analysis and high dimensional
convexity.

Until the present work no techniques were available to approach those
problems nor, more to the point, to determine what the correct questions
were. In this paper we elucidate
this circle of ideas and, in particular, we answer
Problems 1-3 from \cite{M2} negatively. Moreover,
we ascertain the following phenomenon: passing to large subspaces
or quotients can not, in general, erase $k$-dimensional features of a
space if $k$ is below certain threshold  value, which depends on the
dimension of the initial space and the exact meaning of ``large."
For example, the threshold dimension is (at least) of order
$\sqrt{n}$ in the context of ``proportional" subspaces  or quotients
of $n$-dimensional spaces. ``Impossibility to erase" may mean,
for instance, that {\em every } subspace of our $n$-dimensional space
of dimension $\ge~n/2$ will contain a further $1$-complemented subspace
isometric to a preassigned (but a priori arbitrary)
$k$-dimensional space $W$. In a sense, the original space
is ``saturated" with  copies of $W$.
However, this is really just a sample result; our methods are very
flexible -- one may say, canonical -- and clearly can be used
to treat other similar questions.

We employ probabilistic arguments, the most basic idea of which goes
back to Gluskin \cite{G} (see also \cite{MT1} for the survey of other
results and methods in this direction).
Our technique introduces  several novel twists as,
for example, decoupling otherwise not-so-independent events.
In particular, we exhibit families of spaces (or convex bodies)
which enjoy appropriate properties with probability close to $1$.
However,  among other
similar constructions, ours appears to be a relatively good
candidate for ``derandomization" in the spirit of \cite{KT-J} or
\cite{BS}, i.e., for  coming up with an argument that would
yield an explicit space (resp., a convex body) with the same properties.

The organization of the paper is as follows.  In the next section we
state our principal results and their immediate consequences.  Section
\ref{proof} is devoted to the proof of our main result, a general
saturation theorem which implies an  answer to Problem 1 from \cite{M2}.
Section \ref{further} contains the proof of (a version of) Theorem
\ref{hmmp}, relevant to Problems 2 and 3 from \cite{M2}.
Finally, in the Appendix we present a proof of Lemma \ref{shrinking},
a result describing known phenomenon, for which we provide a compact
argument yielding optimal or near optimal constants.

We use the standard notation of convexity and geometric functional
analysis  as can be found, e.g., in \cite{MS1}, \cite{P} or \cite{T}.
The following useful jargon is possibly not familiar to a general
mathematical reader.  A normed space $X$ is completely described by its
unit ball $K=B_X$ or its norm $\|\cdot\|_X$ and so we shall tend to
identify these three objects.  In particular, we will write
$\|\cdot\|_K$ for the Minkowski functional defined by a
centrally symmetric  convex body $K \subset \R^n$ and denote the resulting
normed space  by $(\R^n, \|\cdot\|_K )$ or just $(\R^n, K )$. The
standard Euclidean norm on $\R^n$ will be always denoted by $|\cdot|$.
({\em Attention: } the same notation may mean elsewhere cardinality of a
set and, of course, the absolute value of a scalar.) We will write
$B_2^n$ for the unit  ball in $\ell_2^n$ and, similarly but less
frequently, $B_p^n$  for the unit  ball in $\ell_p^n$,
$1 \le p \le \infty$.

\smallskip\noindent {\small {{\em Acknowledgement } Most of this
    research was performed while the second named author visited
    Universit\'e Marne-la-Vall\'ee and Universit\'e Paris 6 in the
    spring of 2002 and in the spring of 2003.  She wishes to thank
    these institutions for their support and hospitality.} }

\section{Description of  results}
\label{results}
Our first result is a general saturation theorem.
\begin{thm}
\label{hmm}
Let $n $ and $m_0$ be positive integers with $\sqrt{n \log n}
\le m_0 \le~n$.  Then, for every finite
dimensional normed space $W$ with
$$
\dim  W \le c_1 \min\bigl\{ m_0/\sqrt n, m_0^2/(n \log n)\bigr\}
$$
(where $c_1 >0$ is a universal numerical constant) there exists an
$n$-dimensional
normed space $X$ such that  every
quotient $\tilde{X}$ of $X$ with $\dim \tilde{X} \ge m_0$ contains
a 1-complemented subspace isometric to $W$.
\end{thm}

\noindent The theorem has bearing on Problem 1 from \cite{M2}:

\smallskip \noindent {\em Does every }$n${\em -dimensional normed space
admit a quotient of dimension } $\ge~n/2$ {\em whose cotype} 2
{\em constant is  bounded by a universal numerical constant? }

\smallskip
Let us start with several comments concerning the hypotheses
on $k :={\rm dim }\, W $ and $m_0$ included in  the
statement of Theorem~\ref{hmm}.
If, say,  $m_0 \approx n/2$, then
$k$ of order $\sqrt{n}$ is allowed.  Nontrivial
(i.e., large) values of $k$ are obtained whenever $m_0 \gg \sqrt{n
\log{n}}$; we included the condition $\sqrt{n \log n} \le m_0$ to
indicate for which values of the parameters the assertion of the Theorem
is meaningful. Furthermore, for each ``large" quotient
$\tilde{X}$ of $X$, we can deduce quantitative information on
$C_q(\tilde{X})$, the cotype $q$ constant of $\tilde{X}$.
If,  for example, $W = \ell_\infty^k$,
then, for any corresponding $\tilde{X}$, $C_2(\tilde{X})$ is at least
$\sqrt{k}$;  in particular, if
$m_0$ is ``proportional" to $n$, then
$C_2(\tilde{X})$ is at least of order $\root 4\of{n}$.  Similarly,
$C_q(\tilde{X})$ is at least of order $n^{1/2q}$
for any finite $q$.  The problem from \cite{M2} stated above is thus
answered  negatively in a very strong sense.

\smallskip

The dual $X^*$ of the space from Theorem~\ref{hmm} has the property
that, under the same assumptions on $m$ and $k$, every $m$ dimensional
{\em subspace} of $X^*$ contains a (1-complemented) subspace isometric
to $W^* =:V$. If we choose, again, $W = \ell_\infty^k$, all ``large"
subspaces of $X^*$ contain isometrically $\ell_1^k$ and hence $X^*$
comes close to being a counter\-example to Problems 2 and 3 from
\cite{M2}.  Roughly speaking, those problems asked whether every space
of nontrivial cotype $q < \infty$ contains a proportional subspace of
type $2$, or even just $K$-convex. This is well known to be true if
$q=2$ due to presence of nearly Euclidean subspaces.  Although our
space $X^*$ does not {\em a priori } have any cotype property, a
relatively simple modification of our technique allows to correct this
deficiency for sufficiently large $q$, which results in the following
theorem.

\begin{thm}
\label{hmmp}
There exist $q_0 \in [2, \infty)$ such that, for any $q \in (q_0,
\infty)$  there are constants $\alpha = \alpha_q \in  (0,1)$
and $c=c_q >0$ so that the following holds:
for any positive integers $n$ and $m_0$ with
$n^{\alpha} \le m_0 \le n$ and for any normed space   $V$  with
$$
\dim  V \le c_q \ m_0/n^{\alpha}
$$
there exists an
$n$-dimensional nor\-med space $Y$ whose cotype
$q$ constant is bounded by a function of $q$ and $C_q(V)$ and such that
every subspace $\tilde{Y}$ of $Y$ with $\dim \tilde{Y} \ge m_0$
contains a $1$-complemented subspace isometric to $V$.
\end{thm}

Again,  if
-- in the notation of Theorem \ref{hmmp} --
$m_0$ is ``proportional'' to $n$ and $V = \ell_1^k$ is of the maximal
dimension that is allowed (for fixed $q$),  then the type 2 constant of
any corresponding subspace
$\tilde{Y}$ of $Y$ from the Theorem is at least of order
$n^{(1-\alpha_q)/2}$ (and analogously for any nontrivial
type $p > 1$).  The
$K$-convexity constant of any such $\tilde{Y}$ is at least of order
$\sqrt{\log{n}}$ .  Clearly,  this answers in the negative
Problems 2 and 3 from \cite{M2}.

We also remark that choosing $V = \ell_p^k$ (for some $1 < p <2$) in
Theorem \ref{hmmp} leads to a space $Y$ whose type $p$ and cotype
$q$ constants are bounded by numerical constants and such that,
for every $m$-dimensional subspace $\tilde{Y}$ of $Y$ and every
$p < p_1 < 2$, the type $p_1$ constant of $\tilde{Y}$ is at least
$k^{1/p - 1/p_1}$.
If $m_0$ is
``proportional" to $n$, the type $p_1$ constant of $\tilde{Y}$ is at
least of order $n^{(1-\alpha_q)(1/p -   1/p_1)}$, in particular
it tends to $+\infty$  as $n \to \infty$.

\smallskip

In Section \ref{further} we will sketch a proof of Theorem \ref{hmmp}
which gives $q_0=4$.  The Theorem holds in fact even with $q_0=2$, but
the proof, while still based on the same main idea, is much more
subtle and will be presented elsewhere.

\medskip It has been realized in the last few years (cf. \cite{MS2})
that many {\em local } phenomena (i.e., referring to subspaces or
quotients of a normed space) have {\em global } analogues, expressed
in terms of the entire space.  For example, a ``proportional" quotient
of a normed space corresponds to the Minkowski sum of several
rotations of its unit ball. Dually, a ``proportional" subspace
corresponds to the intersection of several rotations.  (Such results
were already implicit, e.g., in \cite{K}.)  In the present paper we
state (without proof) the following sample theorem about the Minkowski
sum of two rotations of a unit ball.

\begin{thm}
\label{hmm-global}
There exists a constant $c_2>0$ such that for any
positive integers $n, k$ satisfying $k \le c_2 n^{1/4}$  and
for any $k$-dimensional normed space $W$, there
exists a normed space $X= (\R^n, K)$ such that, for
any rotation  $\rho: \R^n \to \R^n$, the normed space $(\R^n, K
+ \rho(K))$ contains a $3$-complemented subspace $3$-isomorphic to $W$.
\end{thm}
Theorem \ref{hmm-global} answers a query directed to the
authors by V. Milman. Its proof will be presented elsewhere.

It may be of interest that saturations similar to
those described in Theorems \ref{hmm}-\ref{hmm-global} (and in
Proposition \ref{hmmtech} in the next section) may be
implemented simultaneously for entire families of spaces $W$
whose dimension is uniformly bounded.
For example, if $\mathcal{W}$ is such a family and $|\mathcal{W}|$
is bounded  by a power of $n$ (the dimension of the space that is being
saturated), the strengthening is straightforward: if every
$W \in \mathcal{W}$ is of dimension allowed by a Theorem,
then the space $X$ (or $Y$), whose existence is asserted in the
Theorem, verifies the saturation condition in the assertion
simultaneously for all $W \in \mathcal{W}$. [Only the values of
the numerical constants $c_1, c_q, c_2$ included in the bound on the
dimension $k$ will be affected.]  
In a slightly different direction, given $k \in \N$ and $\ep \in
(0,1)$, it is possible to construct a space of dimension $\le
\exp{(C(\ep) k)}$ whose every ``proportional" subspace (or quotient)
contains a $(1+\ep)$-complemented $(1+\ep)$- isomorphic copy of {\em
  every } $k$-dimensional space.
(The $\exp{(C(\ep) k)}$ dimension estimate is  already optimal for a
super-space  that is universal for all $k$-dimensional spaces.)
We provide a sketch of the argument in a remark at the end of section
\ref{proof}.

\section{The main construction}
\label{proof}

The statement of the main Theorem \ref{hmm} is a particular case
(namely, with $N \approx n^{3/2} \log{n}$) of the following technical
proposition. Recall that for a normed space $W$, we
denote  by $\ell_1^N (W)$  the $\ell_1$-sum of $N$ copies of $W$.

\begin{prop}
       \label{hmmtech}
Let $n $ and $N$
be positive integers with
$ n \log n \le N  \le n^{3/2} \log{n} $.
Let $m_0$ be a positive integer with
$\max\{\sqrt {n \log n}, n^2\log n/N \}\le m_0 \le n$.
Then, for every finite dimensional normed space $W$ with
$$
k :=\dim  W \le c_1  \min \left\{ \frac{m_0}{\sqrt n},
\frac{m_0^2}{n \log n},
\frac{m_0 N}{n^2 \log n}\right\}
$$
(where $c_1 >0$ is an appropriate universal constant) there exists an
$n$-dimen\-sio\-nal normed space $X$ such that, for any $m_0 \le m \le n$,
every $m$-dimensional quotient $\tilde{X}$ of $X$ contains a
1-complemented subspace isometric to $W$. Moreover,
$X$ can be taken as a quotient of $Z = \ell_1^N (W)$.
\end{prop}

\proof
Let $1  \le k \le m \le n \le kN$ be positive integers.
More restrictions  will be added on these parameters
as we proceed.  Notice that choosing the constant $c_1$ small
makes the assertion vacuously satisfied for small values of $m_0$,
and so we may and shall assume that $m_0, n$ and $N$ are large.

Let $W$ be a $k$-dimensional normed space. Identify
$W$ with $\R^k$ in such a way that the Euclidean
ball $B_2^k$ is the ellipsoid of minimal volume containing the
unit ball $B_W$ of $W$. It follows in particular that
$\frac1{\sqrt k} B_2^k \subset B_W \subset  B_2^k$.
Further, this allows to identify  $Z = \ell_1^N (W)$ with $\R^{Nk}$.

Let $G = G(\omega)$
be a $n \times Nk$ random matrix (defined on some underlying
probability space $(\Omega, \pp)$) with independent
$N(0, 1/n)$-distributed Gaussian entries and set
\begin{equation}
\label{def_body}
K= B_{X(\omega)}  := G (\omega) (B_Z) \subset \R^n.
\end{equation}
The random normed space $X = X(\omega)$ is sometimes referred to
as a random (Gaussian) quotient of $Z$, with $G (\omega)$ the
corresponding quotient map.
We shall show that,
for appropriate choices of the parameters, $X(\omega)$
satisfies  the assertion of Proposition \ref{hmmtech} with probability
close to $1$. [The normalization of $G$ is not important;
here we choose it so that the  radius of  the Euclidean ball
circumscribed on  $K$ is typically comparable to 1.]
As usual in such arguments, we will follow the scheme
first employed in \cite{G} which consists of three steps.

\smallskip
\noindent {\em I.\ } For a fixed quotient map $Q : \R^n \ra \R^m$,
the assertion of the theorem holds outside of a small exceptional set
$\Omega^Q$.
\newline {\em II.} The assertion of the theorem is
``essentially stable" under small perturbations of the quotient map,
and so it is enough to verify it for an appropriate net in the set
of all quotient maps.
\newline {\em III.\ } There exists a net 
$\mathcal{Q}$ 
which works in step II such that the set
$\bigcup_{Q\in \mathcal{Q}} \Omega^{Q} $
has small probability.

\smallskip
\noindent By combining steps I and II we see now that the assertion 
of the theorem holds on the complement of the union from step III and
the theorem follows.

\medskip We start by introducing some basic notation.
Denote by  $F_1, \ldots, F_N$
the $k$-dimensional coordinate subspaces of $\R^{Nk}$ corresponding
to the consecutive copies of $W$ in $Z$.
In particular, $B_Z = \conv (F_j \cap B_Z: j=1, \ldots, N)$.
For $j=1, \ldots, N$, we define subsets of $\R^n$ as follows:
$E_j  :=
{G}(F_j) $,
      $K_j := {G}(F_j \cap B_Z)$ and
\begin{equation} \label{K'}
K_j' := {G}(\spn[F_i: i \ne j] \cap B_Z)
= \conv (K_i: i \ne j) .
\end{equation}
Similarly, for $I \subset \{1, \ldots, N\}$,
we let  $K_{I} := \conv (K_i: i \in I) \subset
\R^n $.

\bigskip
\noindent {\em Step I.\ Analysis of a single quotient map. }  
Since a quotient space is determined up to isometry by the kernel of a
quotient map, it is enough to consider quotient maps which are {\em
  orthogonal } projections.  Let, for the time being, $Q: \R^n \to
\R^m$ be the natural projection on the first $m$ coordinates.  In view
of symmetries of our probabilistic model, all relevant 
features of this special
case will transfer to an arbitrary rank $m$ orthogonal projection.

Let $\tilde{G} = Q G$, i.e.,  $\tilde{G}$ is
the $m \times Nk$ Gaussian matrix obtained by restricting $G$ to
the first $m$ rows.
Let $\kk = Q(K) = \tilde{G}(B_Z)$ and denote
the space $(\R^m, \kk)$ by $\tilde{X}$; the space $\tilde{X}$
is the quotient of $X$ induced by the quotient map $Q$.
We shall use the notation of
$\esp_j$,  $\kk_j$,   $\kk_j'$  and $\kk_{I} $
for the subsets of $\R^m$ defined in the same way
as $E_j$, $K_j$, $K_j'$, $K_I$ above,
using the matrix $\tilde{G}$ in place  of $G$.

\medskip
For any subspace $H\subset \R^m$, we will denote by $P_H$
the orthogonal projection onto $H$.
We shall show that outside of an exceptional set of small measure there
exists
$j \in \{1, \ldots, N\}$ such that $P_{\esp_j} ( \kk_j') \subset
\kk_j$.
Since, for any given  $i$, we {\em always } have $\kk = \conv (\kk_i,
\kk_i')$  and $\kk_i \subset \esp_i$, it follows that, for $j$
as above,
\begin{equation}
P_{\esp_j}( \kk) = \conv (\kk_j, P_{\esp_j}( \kk_j')) = \kk_j
      = \esp_j \cap \kk_j.
\label{koniec}
\end{equation}
As $\kk_j$ is an affine image of the ball  $F_j \cap B_Z$,
which is the ball $B_W$ on coordinates from $F_j$,
we deduce  that $\esp_j $ considered as a subspace of
$\tilde{X}$ is then
isometric to $W$ and, moreover, $1$-complemented.
We note in passing that for that subspace to be just
isometric to $W$ (and not necessarily complemented),
a weaker condition $\esp_j \cap \kk_j' \subset
\kk_j$ suffices. This condition can also be analyzed by the
methods of the present paper
(it requires a generalization of Lemma \ref{combin}),
but no improvement results in the range of parameters
that we are interested in.

\medskip
The precise definition of the exceptional set will be somewhat
technical. We start by introducing more subsets of $\R^m$.
Let, for $j =1, \ldots, N$,
$\dd_j := \tilde{G}(F_j \cap  B_2^{Nk}) $ and
\begin{equation}
       \label{dd}
       \dd_j' := \conv  \left(\dd_i:  i \ne j\right).
\end{equation}
That is, $\dd_j'$ is obtained by replacing  $\kk_i$ with $\dd_i$
in the definition of  $\kk_j'$.  We define similarly $\dd$,  $\dd_{I}$
and, analogously,  the subsets $D_j$,  $D_j'$, $D_{I}$ and $D$ of
$\R^n$.  
So, for example, $D_j := {G}(F_j \cap  B_2^{Nk}) $, and 
$D :=\conv (D_j : 1\le j\le N)$.
Note that since $\frac1{\sqrt k} B_2^k \subset 
         B_W \subset  B_2^k$,
it follows that $\frac 1 {\sqrt{k}}D_j \subset K_j \subset D_j$.
Consequently, analogous inclusions hold for all the corresponding
$K$- and $D$-type sets as they are images (or convex hulls, or images
of convex hulls) of the appropriate $K_j$'s  and $D_j$'s.
In particular,  in order for the inclusion
      $P_{\esp_j} (\kk_j') \subset \kk_j$ to hold it is enough to have
\begin{equation} \label{brutal}
P_{\esp_j}( \dd_j') \subset \frac 1 {\sqrt{k}}\dd_j ,
\end{equation}
and it is this seemingly ``brutal" condition that we shall try to enforce.
As a consequence, our results will not be always sharpest possible;
improvements are possible if more geometric information about the
space $W$ (notably its  Banach-Mazur distance to the Euclidean space) is
available.  On the other  hand,  the construction of the
space $X$ will be,  in a sense, universal: we will obtain quotient
maps $G=G (\omega) : \R^{Nk} \ra \R^n$ which, simultaneously
for all $W$, yield spaces verifying the assertion of Theorem \ref{hmm}
(that is,  the resulting spaces are quotients of $\ell_1^N (W)$ via
{\em the same } quotient map).

\medskip Let us go back to the definition of the exceptional set.
      We start by introducing, for $j \in \{1, \ldots,
N\}$, the ``good" sets
\begin{eqnarray}
\label{Omega'(j)}
\Omega_j' &:= &
\left\{ \omega \in \Omega :
P_{\esp_j} (\dd_j') \subset
\kappa B_2^m \right\}\\
\Omega_{j,0}' &:= & \left\{ \frac12 \sqrt{\frac m{n}}(B_2^m \cap
\esp_j) \subset \dd_j
\subset
2\sqrt{\frac m{n}}(B_2^m \cap \esp_j) \right\}, \ \ \
\label{Omega''(j)}
\end{eqnarray}
where $\kappa \in (0,1)$  will be specified at the end of this proof.
Next we set
\begin{equation} \label{Omega(j)}
\Omega_j := \Omega_j' \cap \Omega_{j,0}'.
\end{equation}
      Now if  $\kappa $, $k$,  $m$ and $n$ satisfy
\begin{equation} \label{cubeball}
\kappa \leq \frac 1
        {\sqrt{k}} \cdot \frac12 \sqrt{\frac
m{n}} \, ,
\end{equation}
then, for $\omega \in \Omega_j$, the inclusion
(\ref{brutal}) holds;
{\em a fortiori, } $P_{\esp_j} (\kk_j') \subset \kk_j$ and so
the argument presented earlier applies. Thus
outside of the  exceptional set
\begin{equation} \label{Omega0}
\Omega^0 := \Omega \setminus \bigcup_{1\le j \le N} \Omega_j
= \bigcap_{1\le j \le N}
\left(  \Omega \setminus \Omega_j \right)
\end{equation}
there exists $j \in \{1, \ldots, N \}$  such that
$\esp_j $, considered as a subspace of
$\tilde{X}$,  is isometric to  $W$ and 1-complemented.

\smallskip

It remains to show that the exceptional set $\Omega^0$ is
appropriately small;
this will be the most involved part of the argument.
For $j, i \in \{1, \ldots, N\}$ and $I \subset \{1, \ldots, N\}$
we set
\begin{equation}
\label{Omega(jA)}
\Omega'_{j,I} := \left\{ \omega \in \Omega :
P_{\esp_j}( \dd_I ) \subset
\kappa B_2^m \right\} , \ \ \ \Omega'_{j,i} := \Omega'_{j,\{i\}}.
\end{equation}
In particular, for any $j \in \{1, \ldots, N\}$,
$\Omega_j' = \bigcap_{i\in \{1,\ldots ,N\} 
         \setminus \{ j\}} \Omega'_{j,i}$
(cf. (\ref{dd}), (\ref{Omega'(j)})) and so,  in
view of  (\ref{Omega(j)}),
$$
\Omega_j= \bigcap_{i\in \{0,1,\ldots ,N\} 
       \setminus \{ j\}} \Omega'_{j,i}
$$
This allows rewriting (\ref{Omega0})  as
\begin{equation}
\label{Omega0(1)}
\Omega^0
= \bigcap_{1\le j \le N} \bigcup_{i\in \{0,1,\ldots ,N\} \setminus \{ j\}}
\left(  \Omega \setminus \Omega'_{j,i} \right)
= \bigcup_{(i_j)}\bigcap_{1\le j \le N} \left(
\Omega \setminus \Omega_{j,i_j}' \right),
\end{equation}
where the second union is extended over all sequences
$(i_1, \ldots ,i_N)$ satisfying 
$i_j\in \{0,1,\ldots ,N\} \setminus \{ j\}$
for $j=1,\ldots,N$. The following lemma enables  us to  handle components
of that union by appropriately grouping them.

\begin{lemma}
\label{combin}
Let $\Lambda = (\lambda_{ij})$ be an $N \times N$ matrix such that
\newline $1^\circ$ all its elements are either 0 or 1
\newline $2^\circ$ each column contains at most one 1
\newline $3^\circ$ the diagonal consists of 0's.
\newline Then there exists $J \subset \{1, \ldots, N\}$ such that
$|J| \ge N/3$  and
$$
i, j \in J \Rightarrow \lambda_{ij}=0
$$
\end{lemma}

\noindent
{\em Proof } At the most fundamental level, this lemma can be deduced
from Tur\'an's theorem on existence of independent sets of vertices in
graphs with few edges, cf. \cite{AS}, p.~82: consider a (directed)
graph $G_\Lambda$ whose adjacency matrix is $\Lambda$, then
$G_\Lambda$ has at most $N$ edges and the assertion means that the set
of vertices $J$ is independent. The lower bound on $|J|$ follows by
counting edges in the corresponding extremal Tur\'an graphs.

More conveniently, the formulation we use here can be derived from the
result of K.~Ball on suppression of matrices presented and proved in
\cite{BT}.  By Theorem 1.3 in \cite{BT} applied to $\Lambda$, there
exists a subset $J \subset \{1, \ldots, N\}$ with $|J| \ge N/3$ such
that $\sum_{i \in J} \lambda_{i j} < 1$ for $j \in J$, which is just a
restatement of the condition in the assertion of the Lemma.

Finally, we present here a simple self-contained argument which yields
a variant of the Lemma with $N/3$ replaced by $N/4$ in the estimate
for $|J|$ (which would be sufficient for our purposes).  Let $I := \{i
\ : \ \sum_{j=1}^N \lambda_{ij} \le 1\}$. If $|I^c| \ge N/4$, then
a simple counting argument shows that at
least $N/4$ rows of $\Lambda$ consist only of 0's which clearly yields
the assertion.  If, on the other hand, $|I|> 3N/4$, consider a maximal
subset $J \subset I$ for which the condition on entries given in the
assertion holds. By maximality, for every $i \in I \setminus J$ one of
the elements of $\{\lambda_{ij} \ : \ j \in J \} \cup \{\lambda_{ji} \
: \ j \in J \}$ equals 1 (note that $\lambda_{ii} = 0$ by hypothesis)
and so there at least $|I \setminus J| > 3N/4 - |J|$ such elements.
On the other hand, those elements are contained in the union of $|J|$
columns and $|J|$ rows each of which contains at most one nonzero
element, which implies $2|J| > 3N/4 - |J|$ or $|J| > N/4$, as required.
\qed

We now return to our main argument. Let
$(i_1, \ldots,i_N)$ be any sequence
satisfying $i_j\in \{0,1,\ldots ,N\} \setminus \{ j\}$,
      $j=1,\ldots,N$, and corresponding to a component in the
second union in (\ref{Omega0(1)}). Define a matrix
$\Lambda=(\lambda_{ij})_{i,j=1}^N$  by 
$\lambda_{ij} =1$ if $i = i_j$ for some
$j \in \{1,\ldots,N\}$ and $\lambda_{ij} =0$ otherwise.
Then $\Lambda$ satisfies the assumptions of  Lemma \ref{combin} and
let $J$ be the resulting set of indices
such that $j \in J \Rightarrow i_j \not\in J$ and
that $|J| = \ell := \lceil N/3 \rceil$.
It now follows directly from the definitions (\ref{Omega(jA)}) that
$$\bigcap_{1\le j \le N} \left(
\Omega \setminus \Omega_{j,i_j}' \right)
\subset \bigcap_{ j \in J} \left(
\Omega \setminus \Omega_{j,i_j}' \right)
\subset \bigcap_{ j \in J} \left( \Omega \setminus
(\Omega'_{j,J^c} \cap \Omega_{j,0}')\right) =: \Omega_J .$$
(Here and in what follows the complement  $\cdot^c$ is meant
with respect to  $\{1,\ldots,N\}$, the index $0$ playing a
special role.)
In combination with
(\ref{Omega0(1)}), the above inclusions show  that
\begin{equation}
\label{Omega0cover}
\Omega^0 \subset  \bigcup_{J \in \mathcal{J}} \Omega_J,
\end{equation}
where $\mathcal{J} := \left\{J \subset \{1,\ldots,N\} \
: \ |J| = \ell\right\}$.
Our next objective will be to estimate $\pp(\Omega _J)$ for a fixed
$J \in \mathcal{J}$.
There is no harm in assuming that $J = \{1, \ldots, \ell\}$.
To keep the notation more compact, for $1 \le j \le N$  we let
\begin{equation}
     \label{def-Ej}
\mathcal{E}_j := \Omega \setminus
(\Omega'_{j,J^c} \cap \Omega_{j,0}')
      = (\Omega \setminus \Omega'_{j,J^c})
\cup (\Omega \setminus \Omega_{j,0}').
\end{equation}
In particular, $\Omega_J = \bigcap_{ j \in J} \mathcal{E}_j $.
Let us now make the key observation -- which is apparent from the
definitions (\ref{Omega''(j)}) and (\ref{Omega(jA)}) -- that the events
$\mathcal{E}_j $,  for  $1 \le j \le \ell$,  are
{\em conditionally } independent  with respect to  $D_{J^c}$:
once $D_{J^c}$  is fixed,
and hence $\dd_{J^c}$ is fixed as well,
each $\mathcal{E}_j $  depends only  on the restriction $G_{|F_j}$.
In fact, the ensemble
$\{G_{|F_j} : 1 \le j \le \ell\} \cup \{D_{J^c}\}$ is
independent since its  distinct elements depend on disjoint sets of
columns of $G$,
and the columns themselves are independent.
This and the symmetry in the indices $j \in \{ 1, \ldots,\ell\}$
imply
\begin{equation} \label{decoupling}
\pp(\Omega_J \mid D_{J^c} )
=\pp(\bigcap_{j \in J} (\mathcal{E}_j
\mid D_{J^c} )
= \prod_{j \in J} \pp(\mathcal{E}_j
      \mid D_{J^c} )
=\Bigl(\pp( \mathcal{E}_1  \mid D_{J^c} )\Bigr)^\ell.
\end{equation}
As suggested above, somewhat informally we may think of the above
conditional probabilities as being calculated
``for fixed $D_{J^c}$"
(and hence ``for fixed $\dd_{J^c}$").

\medskip 
In the sequel we shall need
the following auxiliary well-known result on rectangular Gaussian
matrices (see, e.g., \cite{DS}, Theorem 2.13).
We recall that for a $m \times k$ matrix $B$ with $1\le k \le m$,
its singular numbers $s_j(B)$ are the eigenvalues of
$(B^*B)^{1/2}$ arranged in the non-increasing order.
In particular, for all $x \in \R^k$,
$$
s_k(B) |x| \le |Bx| \le s_1(B) |x| .
$$
\begin{lemma}
\label{rectangular}
Let $k, m$ be integers with $1\le k \le m$ and let
$A = (a_{ij})$ be an
$m \times k$ random matrix with independent
$N(0, \sigma^2)$-distributed Gaussian entries and let $t>0$. Then
$$
\pp\left(s_1(A) > (\sqrt{m} + \sqrt{k})\sigma + t \right)
\le e^{-t^2/2\sigma^2} ,
$$ $$
\pp\left(s_k(A) < (\sqrt{m} - \sqrt{k})\sigma - t \right)
\le e^{-t^2/2\sigma^2}.
$$
\end{lemma}

Returning to our main argument, we note that, by the definition
of  $\mathcal{E}_1$,
$\pp( \mathcal{E}_1  \mid D_{J^c} ) \le
\pp(\Omega \setminus \Omega'_{1,J^c} \mid D_{J^c})
+ \pp(\Omega \setminus \Omega_{1,0}'  \mid D_{J^c} )$.
Next, since $\Omega_{1,0}'$ is independent of $D_{J^c}$, the second
term equals just $1 - \pp( \Omega_{1,0}')$. Further, the
condition from (\ref{Omega''(j)}) defining $\Omega_{1,0}'$
can be restated as
$$
\frac12 \sqrt{\frac m{n}} \le s_k(\tilde{G}_{|F_1})
\le s_1(\tilde{G}_{|F_1}) \le 2\sqrt{\frac m{n}}
$$
and hence the probabilities involved can be estimated by using
Lemma \ref{rectangular} with $\sigma=1/\sqrt{n}$ and
appropriate values of $t$. In particular, if $m \ge 16k$, we~get
\begin{eqnarray}\label{bigcube}
\pp(\Omega \setminus \Omega_{1,0}'  \mid D_{J^c} )
&\le&\pp\left(s_k(\tilde{G}_{|F_1}) < \frac12 \sqrt{\frac m{n}}\right)
+ \pp \left(s_1(\tilde{G}_{|F_1}) > 2\sqrt{\frac m{n}} \right)
\nonumber \\
&\le & e^{-m/32} + e^{-9m/32}.
\end{eqnarray}

To estimate the term $\pp(\Omega \setminus \Omega'_{1,J^c} 
\mid D_{J^c})$,
we need to introduce an auxiliary exceptional set
\begin{equation} \label{omega1def}
\Omega^1 := \{ \omega : D \not\subset 2 B_2^n \}
\end{equation}
(recall that the set $D$ was defined  after (\ref{dd}) and that 
$D \supset K$).
An argument parallel to the one that led to (\ref{bigcube})
shows that if $k\le n/16$, then
\begin{equation} \label{omega1}
\pp(\Omega^1) \leq N e^{-9n/32}.
\end{equation}
(An alternative argument uses the Chevet-Gordon inequality --
cf. the proof of Lemma \ref{shrinking} given in the Appendix --
to estimate the expected value of the norm
$\|G : \ell_1^{N}(\ell_2^k) \ra \ell_2^n \|$, and then the
Gaussian isoperimetric inequality to majorize the probablility
that that norm is $>2$, the condition equivalent to that
in the definition of $\Omega^1$.)
For future reference, we  emphasize that
$\Omega^1$ does not depend on $J$  nor
on the projection $Q$.

\medskip
Next we define
$$
\Omega' := \{ \omega : D_{J^c} \not\subset 2 B_2^n \}.
$$
Clearly $\Omega' \subset \Omega^1$.  As we are dealing
with conditional probabilities  $\pp(\cdot \mid~D_{J^c}~)$,
we may -- by the remark following (\ref{decoupling})
-- consider $D_{J^c}$ fixed.  Furthermore,
for the time being we shall restrict our attention to
$\omega \not\in \Omega'$, i.e., to sets
$D_{J^c} \subset 2 B_2^n$,  which implies
$\dd_{J^c} \subset 2 B_2^m$.  Given that  $\Omega'$ is
$D_{J^c}$-measurable, this restriction will not interfere with the
conditional independence of $\mathcal{E}_j$ for $j \in J$.

\medskip

The exceptional sets such as  $\Omega'$ or $\Omega^1$ involve
conditions on diameters of random images of sets.
It is well known that such quantities are related to
the functional $M^*(\cdot)$, defined for a
set $S \subset \R^s$ via
\begin{equation}
     \label{Mstar}
M^*(S) :=
\int_{S^{s-1}}{\sup_{y \in S}
{\langle x, y \rangle}} dx ,
\end{equation}
where the integration is performed with respect to the normalized
Lebesgue measure on $S^{s-1}$ (this is $1/2$ of what  geometers
call the mean  width of $S$). It is elementary to verify that
for $\omega \not\in \Omega^1$
(hence, {\em a fortiori}, for $\omega \not\in \Omega'$)
and for any $j = 1, \ldots, N$,
$$
M^*(\dd_j) \le 2 M^*(\esp_j \cap B_2^m)
= 2 \int_{S^{m-1}} |Px| dx
\le 2\sqrt{k/m},
$$
where  $P = P_{\esp_j}$ is the orthogonal projection on
$\esp_j$ (or, by rotational
invariance,  {\em any} orthogonal projection of rank $k$).
It follows then via known arguments (most readily,  by passing to
Gaussian averages, cf. \cite{LT}, (3.6), or Lemma 8.1 in \cite{MS})
that, again for  $\omega \not\in \Omega^1$,
\begin{equation}
       \label{tled}
M^*(\dd_{J^c}) \le C_0 \sqrt{\frac{\log{N}}m} + \max_{j \in J^c}
{M^*(\dd_j)} \le (C_0 +2)\sqrt{\frac{\max\{k, \log N\}}{m}},
\end{equation}
where $C_0\ge 1$ is a  universal constant.

\medskip
Let $d, m$ be integers with  $1 \le d \leq m$ and
let $G_{m,d}$ be the Grassmann manifold of $d$-dimensional
subspaces of $\R^m$ endowed with the normalized Haar measure.
The following lemma describes the behavior of the diameter
of a random rank $d$ projection of a subset of  $\R^m$.

\begin{lemma}
\label{shrinking}
Let  $a >0$ and let $S \subset \R^m$ verify \ $S \subset a B_2^m$.
Then, for any  $t > 0$, the set
$\left\{ H \in G_{m,d} : P_H (S) \subset
\left(a \sqrt{d/m} + M^*(S) + t \right) B_2^m \right\}$
has measure $\geq 1- \exp(-t^2m/2a^2 + 1)$.
\end{lemma}
The phenomenon
discussed in the Lemma is quite well known, at least if
we do not care about the specific values of numerical constants (which
are not essential for our argument) and precise estimates on
probabilities.  It is sometimes refered to as the ``standard
shrinking'' of the diameter of a set, and it is implicit, for example,
in probabilistic proofs of the Dvoretzky theorem, see \cite{M-dv}, 
\cite{MS1}.
For future reference, we provide in the Appendix a compact
proof which yields the formulation given above.

\medskip Returning to our main argument,
we now use the  estimate from Lemma~\ref{shrinking} for
$S = \dd_{J^c}$, $d = k$ and $t = \kappa/2$.
Since $\omega \not\in \Omega'$, we may take $a =2$
to  deduce  that
the set
$\{ H \in G_{m,k} : P_H (\dd_{J^c}) \subset
(2 \sqrt {k/m} + M^*(\dd_{J^c}) + \kappa/2)
B_2^m \}$ has measure larger than or equal to
$ 1- \exp(-\kappa^2 m / 32 +1)$.
Taking into account (\ref{tled}) and requiring that,
in addition to (\ref{cubeball}),
  $\kappa $ verifies the condition
\begin{equation}
      \label{k1}
C' \sqrt{\max\{k,\log{N}\}/m} \le \kappa
\end{equation}
with
\begin{equation}
       \label{const_k1}
       C' \ge 2 (C_0 +4),
\end{equation}
we see that
$$
2 \sqrt {k/m} + M^*(\dd_{J^c}) + \kappa/2 \le \kappa,
$$
and thus the measure of the set
$\{ H \in G_{m,k} : P_H (\dd_{J^c}) \subset
\kappa B_2^m \}$ is also larger than or equal to
$ 1- \exp(-\kappa^2 m / 32 +1)$.
Since, by the rotational invariance of the Gaussian
measure, the distribution of $\esp_1$ on $G_{m,k}$ is uniform,
we finally get  that
\begin{eqnarray}
      \label{smallball}
\pp\left(\Omega'_{1,J^c}  \mid D_{J^c} \right)
&=&\pp\left(P_{\esp_1}( \dd_{J^c}) \subset
\kappa B_2^m  \mid D_{J^c}
      \right) \nonumber \\
&\geq &1- \exp (-\kappa^2 m / 32 +1).
\end{eqnarray}
The above  estimate
is valid outside of $\Omega'$, i.e.,
for $D_{J^c} \subset 2 B_2^n$. Thus, recalling
(\ref{def-Ej}) and combining
(\ref{smallball}) and (\ref{bigcube})
    we infer that,
under the same restriction,
$$
\pp(\mathcal{E}_1  \mid D_{J^c} )
\le e^{-\kappa^2 m/32 +1}+ e^{-m/32} + e^{-9m/32}.
$$
Plugging the above into (\ref{decoupling}) and recalling that,
     by (\ref{cubeball}), $\kappa \le 1/2$,  we obtain
$$
\pp(\Omega_J \mid D_{J^c} ) \leq (2e\,e^{-\kappa^2 m/32})^ \ell =
(2e)^\ell e^{-\kappa^2 m \ell/32},
$$
again,  on the set $\Omega \setminus \Omega'$ which is
      $D_{J^c}$-measurable.
Consequently, averaging over $\Omega \setminus \Omega'$,
$$
\pp(\Omega_J \mid \Omega \setminus \Omega')
      \leq (2e)^\ell e^{-\kappa^2 m \ell/32}.
$$
Since (as was noted following the definition of $\Omega'$)
$\Omega' \subset \Omega^1$, our argument shows that
\begin{equation}
 \label{OmegaJ}
\pp(\Omega_J \setminus \Omega^1) \leq \pp(\Omega_J \setminus
\Omega') = \pp(\Omega_J \mid \Omega \setminus \Omega')
\pp(\Omega \setminus \Omega') \leq  (2e)^\ell e^{- \kappa ^2 m \ell/32}.
\end{equation}

Up to now the set $J$ was fixed,  now is the time to allow it
to vary over $\mathcal{J}$. Since
$|\mathcal{J}| = {{N}\choose{\ell}}$,  it follows that
$$
\pp\left(\bigcup_{J \in \mathcal{J}} \Omega_J \setminus \Omega^1
\right) \leq {{N}\choose{\ell}}  (2e)^\ell  e^{- \kappa^2 m \ell/32},
$$
and so,  by (\ref{Omega0cover}),
\begin{equation} \label{singleop}
\pp(\Omega^0) \leq \pp(\Omega^1) + \pp\left(\bigcup_{J \in
\mathcal{J}} \Omega_J \setminus \Omega^1
\right) \leq N e^{-9n/32} + {{N}\choose{\ell}} (2e)^\ell
               e^{-\kappa^2 m \ell/32}.
\end{equation}

Later on we will make the two terms on the right hand side small, by
appropriate assumptions on $k$ and appropriate choices of $\kappa$
and $N$. Then, recalling the definition of $\Omega^0$, we will
deduce that $\tilde{X}$, which is the quotient of the random space $X$
via the map $Q$, with large probability contains a 1-complemented
subspace isometric
to $W$. However, our ultimate goal is to show that this is true for
{\em all} $m$-dimensional quotients of $X$.  As indicated earlier the
strategy for such an argument depends on combining a sharp probability
estimate for a single $Q$ (obtained in Step I above) with perturbation
and discretization arguments involving an appropriate net in the set
of all $Q$'s (shown in Steps II and III below).

Finally, we note that
the condition $m \ge 16 k$ that was used in deriving
estimates  (\ref{bigcube}) and  (\ref{omega1}) is implied by
  (\ref{k1}) combined with  (\ref{const_k1}), and so
in what follows we need to ensure only the latter two
conditions.

\bigskip
\noindent
{\em Step II.\ The perturbation  argument }
We first generalize the definition
(\ref{Omega0}) of the exceptional
set.  Let $Q$ be {\em any }  orthogonal projection of rank\,$m$
(which will remain fixed through the end of the present step).
Let us denote by $\Omega^Q$ the set given by
formally the same formulae as in (\ref{Omega0}) by the Gaussian
operator $\tilde{G} = QG$ for this particular $Q$.
By rotational invariance,  all the properties we
derived for $\Omega^0$  hold also for $\Omega^Q$.
The object of Step II will be to show that --
under appropriate hypotheses -- properties just slightly weaker,
but still adequate for our purposes, hold for projections
sufficiently close to $Q$.

\smallskip
Recall that by the definition (\ref{omega1def}) of $\Omega^1$ for
$\omega \not\in \Omega^1 $
we have, for $j \in \{1, \ldots , N\}$,
\begin{equation} \label{global}
D_j' \subset D \subset 2 B_2^n,
\end{equation}

Let now $\omega \not\in \Omega^1 \cup \Omega^Q$
and let $Q'$ be a rank\,$m$ projection  such that
$\|Q-Q'\| \le \delta$, where $\| \cdot \|$ is the operator
norm with respect to $|\cdot |$ and $\delta > 0$ a constant to be
specified later.  We recall that our objective is to show that,
for some $j$, conditions just slightly weaker than those in
(\ref{Omega'(j)}) and (\ref{Omega''(j)}) hold
with $Q$ replaced by $Q'$.  Consequently, one still will be
able to deduce that the quotient of $X$  corresponding to
$Q'$ also contains a 1-complemented subspace isometric to $W$,
namely $Q'E_j$.

\smallskip
Since $\omega \not\in \Omega^Q$, it follows from
(\ref{Omega(j)}) and  (\ref{Omega0}) that
  $\omega \in \Omega_j
= \Omega_j' \cap \Omega_{j,0}'$ for some $j \in \{1, \ldots , N\}$, and
so  the conditions from
(\ref{Omega'(j)}) and (\ref{Omega''(j)}) hold for this particular
$j$ and $Q$. (Note that $Q$ enters implicitly into definitions of
all ``tilde-objects''.)
Let us first analyze
(\ref{Omega''(j)}).  It  asserts that,  for
$x \in F_j$, we have two-sided estimates
\begin{equation}
       \label{net''}
\frac12 \sqrt{\frac m{n}} |x| \le |\tilde{G}x|=
|QGx| \le  2\sqrt{\frac m{n}} |x|.
\end{equation}
We want to show  similar estimates for $Q'$.  We have
$$
      \Big| \ |Q'Gx| - |QGx|\ \Big| \le \|Q'-Q\| \, |Gx|
\le \delta |Gx| .
$$
On the other hand,  $\omega \not\in \Omega^1$ is equivalent
to $\|G_{|F_i}\| \le 2$ for all $1 \le i \le N$ (as an operator
with respect to  the Euclidean norms; see (\ref{omega1def})
and the paragraph
following (\ref{omega1})). Thus, by (\ref{net''}), we  get
$(1/2) \sqrt{m/n}|x| -  2\delta |x|\le |Q'Gx|$,
and an analogous upper estimate, for all   $x \in F_j$.
So if,  say, $\delta \le  (1/8) \sqrt{m/n}$,
we obtain an analogue of (\ref{net''}) -- 
hence of the condition  from    
(\ref{Omega''(j)}) --  for $Q'$, namely
\begin{equation}
  \label{netP''}
\frac14 \sqrt{\frac m{n}} |x| \le
|Q'Gx| \le  \frac 94\sqrt{\frac m{n}} |x|
\end{equation}
for all  $ x \in F_j$.

We now turn to (\ref{Omega'(j)}).  It asserts that for $y \in D_j'$
we have $|P_{\esp_j}(Qy)| \le \kappa$ and, again, we want
to prove a similar estimate for $Q'$.  To make the
notation more compact, set $\bar E_j := Q'G(F_j)$ for $j=1,
\ldots, N$.  We shall show first that
\begin{equation}
       \label{dif_qs}
\|P_{\esp_j} - P_{\bar E_j}\| \le 4 \delta \sqrt{\frac{n}{m}} 
=: \delta_1.
\end{equation}

It is a general and elementary fact that the difference $P_{H} -
P_{H'}$ of two orthogonal projections of the same rank attains its
norm on {\em each} of the ranges of the projections.  (To see this,
consider the Schmidt decomposition of the restriction $P_{H|H'}=
\sum_i \lambda_i \langle \cdot, h_i'\rangle  h_i$ and observe that the
mutually orthogonal subspaces $\spn (h_i, h_i')$ are reducing for
$P_{H} - P_{H'}$, and that $P_{H} - P_{H'}$ is a multiple of an
isometry on each of these subspaces.)  It is thus enough to estimate
the norm of $P_{\esp_j} - P_{\bar E_j}$ restricted to the subspace
$\esp_j$.  If $u = QGx \in \esp_j$, with $x \in F_j$, then
$$
(P_{\esp_j} - P_{\bar E_j})u =
(I- P_{\bar E_j}) ( Q-Q') Gx.
$$
Therefore, using (\ref{net''}), we get, for  $u \in \esp_j$,
$$
|(P_{\esp_j} - P_{\bar E_j})u |\le \|Q - Q'\| \|G_{|F_j}\||x|
\le 2 \delta (2\sqrt{n/m})|u|,
$$
and (\ref{dif_qs}) immediately  follows.

\medskip
Now let $y \in D_j'$; then, by (\ref{global}), 
 $|y|\le 2$ and
\begin{eqnarray*}
|P_{\bar E_j} Q' y|
& \le & |(P_{\bar E_j} - P_{\esp_j} )Q'y |
+ |P_{\esp_j} (Q' - Q) y |  +  |P_{\esp_j} Qy|\\
& \le & 2 \delta_1 + 2 \delta +  \kappa.
\end{eqnarray*}
Thus if, say, $\delta_1 \le \kappa/4$, then -- since clearly $\delta <
\delta_1$ -- it follows that $ |P_{\bar E_j} Q' y| \le 2 \kappa$, and
so the analogue of (\ref{Omega'(j)}) holds for $Q'$ with $2\kappa$
replacing $\kappa$.

\smallskip

We now set $\delta := 1/(8 \sqrt n)$. It is then easy to see that the
above condition for $\delta_1$ holds, just combine (\ref{k1})
and  (\ref{const_k1}); the requirement $\delta \le (1/8)\sqrt{m/n}$,
discussed earlier  in the context of  the analogue of
(\ref{net''}), is then also trivially satisfied. We thus conclude that
if $\omega \not\in \Omega^1 \cup \Omega^Q$,
$\|Q-Q'\| \le \delta $ and
\begin{equation} \label{cubeball4}
2 \kappa
\leq \frac 1{\sqrt{k}} \cdot \frac14 \sqrt{\frac m{n}} ,
\end{equation}
then the quotient of $X$  corresponding to
$Q'$ contains an  isometric and 1-complemented copy of $W$.
We note that (\ref{cubeball4}) is just slightly stronger than
(\ref{cubeball}),  and as easy to satisfy.

\bigskip
\noindent
{\em Step III.\ The discretization: a net in the set of quotients }
The final step is now straightforward. If $\mathcal{Q}$
is any finite family of rank$\, m$ orthogonal
projections,  then
\begin{eqnarray} \label{manyop}
\lefteqn{
\pp\left(\Omega^1 \cup 
\bigcup_{Q\in \mathcal{Q}} \Omega^Q \right)
= \pp\left(\Omega^1 \cup 
\bigcup_{Q\in \mathcal{Q}}
(\Omega^Q \setminus \Omega^1)\right) }\\
&\leq &  N e^{-9n/32} + 
|\mathcal{Q}| \, {{N}\choose{\ell}}
(2e)^\ell
e^{-\kappa^2 m \ell/32} \nonumber
\end{eqnarray}
as in (\ref{singleop}); just note that the exceptional set
$\Omega^1$ 
does  not depend on $Q$.
If $8 \log N \le n$ the first term is
$< e^{-n/8}$.
Now recall that  the set of rank $ m$ orthogonal projections on
$\R^n$, endowed with the distance given by the operator norm, can be
identified with $G_{n,m}$ endowed with the appropriate invariant
metric.  The latter set admits, for any $\delta>0$, a $\delta$-net
$\mathcal{N}$ of cardinality $|\mathcal{N}| \le
(C_2/\delta)^{m(n-m)}$, where $C_2$ is a universal constant  (see,
e.g., \cite{S2}).  For our choice of $\delta = 1/(8 \sqrt n)$,
this does not exceed $e^{mn\log{n}}$,
at least for sufficiently large $n$.  We
plug this estimate into (\ref{manyop}).

Recall that
$\ell = \lceil N/3 \rceil \ge N/3$
and let us  assume that the constant $C'$ in
the definition (\ref{k1}) of $\kappa $ satisfies also $ C'\ge 20$,
so that $\kappa^2 m \ge 400$. Then
$$
{{N}\choose{\ell}}
(2e)^\ell
e^{-\kappa^2 m \ell/32}
\le 2^N \left( 2 e^{1- \kappa^2 m/32} \right)^{\ell} \le
e^{(1+ 4 \log{2}- \kappa^2 m/32)N/3} \le e^{-\kappa^2 m N/128}.
$$
Thus  the 
second  term in (\ref{manyop})
is  less than or equal to
$
e^{mn\log{n} - \kappa^2 m N/128}.
$

\medskip \noindent
In conclusion, if  $k$, $\kappa$,  $m$,  $n$ and $N$ satisfy
\begin{equation}
C' \sqrt{\max \{k, \log{N}\}/m}  \le  \kappa \le 1,
\qquad 2^8\, mn\log{n} \le  \kappa^2 m N,  \qquad 
8  \log N \le n,
\label{warunki}
\end{equation}
where
$$
C' = \max\{20,  2(C_0 + 4)\}
$$
(see (\ref{k1}) and   (\ref{const_k1})),
then the set $\Omega \setminus (\Omega^1 \cup 
\bigcup_{Q\in \mathcal{Q}} \Omega^Q)$ has positive measure (in fact,
very close to 1 for large $n$).
If, additionally, (\ref{cubeball4}) is satisfied, then any $\omega$
from this set induces an $n$ dimensional space $X$ whose {\em all }
$m$-dimensional quotients contain an isometric 1-complemented copy of
$W$.  (Indeed, by the argument above, each such quotient is determined
by a projection within $\delta = 1/(8 \sqrt n)$ of a certain $Q \in
\mathcal{Q}$.)  Then the assertion of Proposition~\ref{hmmtech} holds for
that particular value of $m$.

\smallskip
To ensure that the conditions in (\ref{cubeball4}) and
(\ref{warunki})  are consistent,
it will be most convenient to
let  $\kappa  = ( m/(64\, n k))^{1/2}  $,
so that (\ref{cubeball4}) holds (in fact with equality).
The first inequality,
$ \kappa \ge C' \sqrt{\max \{k, \log{N}\}/m} $, becomes then $k \le c_1'
\min
\{ m/\sqrt n, m^2/(n\log N)\}$, for some numerical constant $ c_1' >0$.
The second   condition in (\ref{warunki}) is equivalent to
$k \le 2^{-14} {m N}/(n^2 \log n)$,
and hence the requirements on $k$ may be summarized as
\begin{equation}
\label{warunki2}
k \le c_1'' \min \left\{ \frac{m}{\sqrt n}, \frac{ m^2}{n\log N},
        {\frac{m N}{n^2 \log n}}\right\},
\end{equation}
where $c_1'' \in (0,1)$ is an appropriate absolute constant.

Only the cases when the fractions on the right hand side of
(\ref{warunki2}) are at least $1$ are of interest.
In particular, $n^2 \log n/ N \le m \le n$,
which gives the lower estimate $N \ge n \log n$.
(Lower bounds on $m$ from the hypotheses of the
Proposition will follow similarly.)
On the other hand, for $N \ge n^{3/2} \log{n}$
the right hand side of (\ref{warunki2}) does not depend on $N$
anymore, hence the interval $[n\log n, n^{3/2} \log{n}]$
does indeed include all interesting values of $N$.
The third   condition in (\ref{warunki}) is  then
easily satisfied once $n$ is large enough.

Since $\log n \sim \log N$,
(\ref{warunki2}) specified to $m=m_0$ reduces to the hypothesis of
Proposition  \ref{hmmtech}, and if it holds for $m=m_0$, it is
necessarily true in the entire range $m_0 \le m \le n$.
It follows that, under our hypothesis, the above construction
can be implemented for each $m$ verifying $m_0 \le m \le n$.  Moreover,
since the estimates on the probabilities of the exceptional sets
corresponding to different values of $m$ are exponential in $-m$ (as
shown above), the sum of the probabilities involved is small.
Consequently, the construction can be implemented {\em simultaneously}
for all such $m$ with the resulting space satisfying the full
assertion of Proposition~\ref{hmmtech} with probability close to $1$.
\qed

\noindent
\rem We now sketch arguments that are behind the comments following
Theorem \ref{hmm-global} and concern simultaneous saturations. They
are based on the following observation: the construction remains
virtually unchanged if we add ``not too many" additional arbitrary
factors of dimension not exceeding $k$ to -- for example -- the
$\ell_1$ sum defining $Z$; see the beginning of the proof of
Proposition \ref{hmmtech}. This is because the number of such factors
affects the estimates rather lightly; in the context of Proposition
\ref{hmmtech}, its only effect on the argument is via $N$ in the
second expression  on the right hand side 
of formula (\ref{warunki2}).  Accordingly, if we
replace $Z=\ell_1^N(W)$ by $\left(\bigoplus_{W \in
  \mathcal{W}}\ell_1^N(W)\right)_1$ (a direct sum in the
$\ell_1$-sense), 
the binomial coefficient $N\choose{\ell}$ in 
(\ref{manyop}) and the argument following it needs to be replaced by 
${{N |\mathcal{W}|}\choose{\ell}} < (e N |\mathcal{W}|/\ell)^\ell  $.
The  corresponding upper bound on $k$ becomes then $c
m^2/(n \log(N |\mathcal{W}|))$. This implies the first comment.  
For the second comment it is enough to invoke the known
fact that the family of $k$-dimensional normed spaces admits a
$(1+\ep)$-net (in the Banach-Mazur distance) whose cardinality
is at most $\exp(\exp(C(\ep)k))$. [This estimate and its optimality
are quite straightforward to derive employing, in particular, the
methodology of \cite{G}, but in fact they have been shown earlier, see
\cite{br}; we thank the referee for providing us with this reference.]
Finally, the optimality of the $\exp(C(\ep)k)$ dimension estimate
for a space containing a $(1+ \ep)$-isomorphic copy of every
$k$-dimensional normed space follows from the known estimates on
nets of Grassmann manifolds, already used in Step III earlier in this
section, and from the optimality of the $\exp(\exp(C(\ep)k))$ estimate
mentioned above.

\section{Subspaces of spaces with finite cotype}
\label{further}

As mentioned in Section \ref{results}, we shall prove
Theorem~\ref{hmmp} for $q >4$ only, in which case it will be an
immediate consequence of the following technical statement.

\begin{prop}
     \label{hmmptech}
Let $q \in (4, \infty)$ and set $\bt := (q+2)/(2q-2)
\in (1/2, 1)$ and $\gamma := (q+1)/(2q-2)$.
Let $n$ and $m_0$ be positive integers with $n^\bt
(\log{n})^{\gamma} \le m_0 \le n$.
Let $V$ be any normed space with
$$
\dim  V \le c_0 \min\Bigl\{
\frac{m_0} {n^{\bt} (\log{n})^{\gamma - 1/2}}, \frac {m_0^2}{ n^{2 \bt}
(\log n)^{2 \gamma}}  \Bigr\}
$$
(where $c_0 \in (0,1)$ an appropriate numerical constant). Then there
exists an
$n$-dimensional nor\-med space
$Y$ whose cotype
$q$ constant is bounded by a function of $q$ and $C_q(V)$ and such that,
for any \ $m_0 \le m
\le n$, every $m$-dimensional subspace $\tilde{Y}$ of $Y$ contains a
$1$-complemented subspace isometric to $V$.
\end{prop}

In fact
we know how to show a version of Proposition \ref{hmmptech} under the
hypothesis $\dim V \le {c_q \ m_0}/(n^{1- \eta} (\log{n})^{
(1-2\eta)/3})$, where $\eta := (q-2)/(2q+2)$ ($\in (0, 1/2)$) and
$c_q >0$ is a constant depending on $q$ only.  This gives a nontrivial
outcome for any $q > 2$, and thus proves Theorem~\ref{hmmp} in the
full generality. However, the argument is much more subtle and
involved than the one presented below.

\medskip \noindent {\em Proof of Proposition \ref{hmmptech} \ }
The idea behind the argument is as follows. The starting point is a dual
$X^*$ of the space from Theorem \ref{hmm} which, in the present notation,
is a subspace of $\ell_\infty^N(V)$, where, in particular, $N < n^2$.
The trick is to consider the same subspace, but this time endowed with
the norm inherited from $\ell_q^N(V)$ (which gives control of the cotype
property). Since the ratio between the
two norms is at most $N^{1/q}$, and since -- as may be verified --
the margin of error in conditions ensuring the assertion of
Theorem \ref{hmm} was also a (small but fixed) power of $n$,
the rest of the proof carries over to the new setting if $q$ is
large enough.

And here are the main points of the argument.
Fix $q >4$
and let $p = q/(q-1)$ be the conjugate exponent.
Let $1\le k \le m \le n \le kN$ be positive integers
verifying conditions (\ref{warunki}).
We shall specify $N$ (depending on $q$) at the end of the proof
when we shall also derive additional restrictions on $k$ and $m$.
If $V$ is such that $\dim V = k$, we
let $Z_p = \ell_p^N (V^*)$ and let $X_p = X_p(\omega)$ be the
Gaussian quotient of $Z_p$ obtained analogously as in the proof
of Proposition~\ref{hmmtech} with $W= V^*$; denote its unit ball by
$K_p$.  Otherwise, we shall use the same notation as that of the
proof of Proposition~\ref{hmmtech}, and we shall make the same choice
of $\ell$.  The dual spaces $X_p^* = X_p(\omega)^*$ are thus
isometric to  subspaces of $Z_p^* = \ell_q^N (V)$ and so their
cotype $q$ constants are uniformly bounded (depending on $q$ and
the cotype $q$ constant of $V$).  We claim that once all the parameters
are properly chosen then, outside of a small exceptional set,
$Y=X_p(\omega)^*$ satisfies also the (remaining) condition of the
assertion of Proposition \ref{hmmptech}
involving the subspaces isometric to $V$.
That condition can be restated as follows:
every quotient of $X_p(\omega)$ of dimension $m \ge m_0$ contains a
$1$-complemented subspace isometric to $V^*$  for values
of $k$  described in Proposition \ref{hmmptech}.
Thus we have a very similar problem to the one encountered in
Proposition \ref{hmmtech}, and the strategy is to show that simple
modifications of the argument applied then do yield the result asserted
in the present context.

\medskip To that end, observe that $F_j \cap B_{Z_p} = F_j \cap B_{Z}$
for all $j \in \{1, \ldots , N\}$ (recall that $F_1, \ldots, F_N$ are
$k$-dimensional coordinate subspaces of $\R^{Nk}$, defined in the
paragraph preceding (\ref{K'})).
Moreover,
$B_Z \subset B_{Z_p} \subset  N^{1/q} B_Z$. Consequently,
the same inclusions hold for the images of these balls
by $G$ and $\tilde{G}$ and their appropriate sections and projections.
Similarly,  for any $j \in \{1, \ldots , N\}$, the set
$$
\tilde{L}_j' := \tilde{G}(\spn[F_i: i \ne j] \cap B_{Z_p})
$$
(cf. (\ref{K'})) satisfies $\tilde{L}_j'
\subset (N-1)^{1/q}\kk_j'$. Hence (cf. (\ref{Omega'(j)}))
$$
\left\{ \omega \in \Omega :
P_{\esp_j} ( \tilde{L}_j')\subset
(N-1)^{1/q} \kappa B_2^m
\right\}
\supset \Omega_j',
$$
where, as before, $\kappa$ satisfies (\ref{k1}).
Accordingly, for the argument of Proposition~\ref{hmmtech} to carry
over to the present setting (with the same exceptional set!), we need
to require, in place of (\ref{cubeball4}), a stronger inequality
\begin{equation} \label{cubeballq}
(N-1)^{1/q} \cdot 2\kappa
\leq \frac 1{\sqrt{k}} \cdot \frac14 \sqrt{\frac m{n}} .
\end{equation}
It remains to reconcile this condition with (\ref{warunki}),
which provides {\em lower } bounds on $\kappa$ that may
be summarized as
\begin{equation}  \label{warunkiq}
\kappa^2 \ge C_0 \max \left\{ \frac{ n \log n} {N}, \frac {k}{m},
\frac{\log N}{m} \right\},
\end{equation}
where $C_0 > 1 $ is a numerical constant.
[These bounds come from the first two inequalities  in
(\ref{warunki}); as before, the third  inequality will
be automatically satisfied.] We now choose
$\kappa$  to verify (\ref{cubeballq}), say,  $\kappa =
\sqrt{m/(64\, n k)} \,N^{-1/q} $, plug this value into
(\ref{warunkiq})
and solve the obtained inequalities for $k= \dim V$.
This leads to the constraint
\begin{equation} \label{warunkiqk}
\dim V \le c_0'
\min \left\{ \frac {m N^{1-2/q}}{ n^2 \log n}, \frac{m}{\sqrt n N^{1/q}},
\frac{m^2} {n N^{2/q}\log N} \right\},
\end{equation}
where $c_0' \in (0,1)$ is a universal constant. It remains to choose $N$
to optimize this constraint. However, we are dealing with a range of
possible values of $m$ and there is no universal optimal choice.
Still, concentrating on the first two expressions in (\ref{warunkiqk}),
for which the best value of $N$ is of order $(n^{3/2} \log n)^p$, gives a
``near optimal" outcome asserted in Proposition \ref{hmmptech}:
just replace everywhere
$N$ by $(n^{3/2} \log n)^p$, note that then $\log N \sim
\log n$, and that the most restrictive condition comes from the
smallest value of $m$, namely $m=m_0$.
(Ignoring the third expression in (\ref{warunkiqk}) when optimizing
$N$ is justified by the fact that it is larger than the second one
except in a narrow ``logarithmic" range of $m$.)  The lower bounds
$n^\bt (\log{n})^{ \gamma}$ on $m$ and $4$ on $q$ are included to
indicate the values for which the Proposition is of interest.  As is
easily seen, if any of them is violated, then the only allowed value
for $\dim V$ is $0$. By contrast, if $q > 4$ and, say, $m_0$ remains a
fixed proportion of $n$, then the upper bound on $\dim V$ grows
without bound as $n \ra \infty$, as needed for non-trivial
applications to Theorem \ref{hmmp}-like statements. \qed

\section{Appendix}
{\em Proof of Lemma~\ref{shrinking} } Let $A = (a_{ij})$ be
a $d \times m$ random matrix with independent $N(0, 1/m)$-distributed
Gaussian entries. We first show that an estimate very similar to the
assertion of Lemma \ref{shrinking} holds if we replace $P_H$ by
$A$, and then present a rather general argument which allows to
deduce the estimates for $P_H$.
(Other variants of a general argument connecting Gaussian and
orthogonal settings were developed recently in \cite{MT2}.)
We shall consider only the case $m > 1$ (the assertion of the Lemma
holds trivially if $m=1$), and we shall assume that $S$ is closed.

We begin by applying the Chevet-Gordon inequality (see, e.g.,
\cite{Go}, Theorem 5) to obtain
\begin{equation} \label{gordon}
{\mathbb{E}} \max_{x \in S}|A x|  \le
\sqrt{\frac dm} \, a  + M^*(S) ;
\end{equation}
[In fact a slightly stronger inequality holds: $\sqrt{d}$
and $ M^*(S)$ may be replaced by the smaller quantities
$\mathbb{E}|\cdot |$ and $\mathbb{E} {\sup_{y \in S} {\langle
x, y \rangle}}/\sqrt{m}$, respectively. The expectations are
being taken \wrt the standard Gaussian measure on the 
appropriate Euclidean space.]
 Next, since the function $A
\ra \max_{x \in S}|A x|$  is $a$-Lipschitz \wrt the
Hilbert-Schmidt norm, it follows from the Gaussian
isoperimetric inequality (see, e.g., \cite{LT})
 that, for any $t > 0$,
\begin{equation} \label{isoper}
\pp\left( \max_{x \in S}|A x| \ge \sqrt{\frac dm} \, a  + M^*(S) + t
\right) \le 
e^{-t^2m/2a^2}.
\end{equation}
(Notice  that the variance of each entry of $A$ is   $ 1/ m$.)
Denote the  normalized Haar measure on $G_{m,d}$ by $\nu$.
Lemma \ref{shrinking} will easily follow if we show that,
for any $\tau > 0$,
\begin{equation} \label{passing}
\pp \left(A :  \max_{x \in S}|A x| \ge
\tau \right) \ge
e^{-1} \, \nu \left(H \in G_{m,d} :  \max_{x \in S}|P_H x| \ge
\tau \right).
\end{equation}
To that end, we note first that, in the calculation of the measure
on the right hand side above,
we may replace $P_H$ by $R_dU$ and $\nu$ by $\mu$, where
$R_d : \R^m  \ra \R^d$ is the restriction to the first $d$ coordinates,
$U \in O(m)$, and $\mu$ is the Haar measure on $O(m)$;
this follows from the invariance of $\nu$ under the action of $O(m)$.
Similarly, in the calculation of the probability we may replace
$A$ by $(A A^*)^{1/2}R_dU$ and $\pp$ by $\pp \otimes \mu$.
Indeed, by the invariance of the Gaussian measure under
the action of $O(m)$, the radial and the polar parts of $A$ are
independent and the latter must have the same distribution as $R_dU$.
In other words, the inequality above can be rewritten as
\begin{equation} \label{polar}
\pp \otimes \mu \left(\max_{x \in S}|(A A^*)^{1/2}R_dU x| \ge
\tau \right) \ge
e^{-1} \, \mu \left(\max_{x \in S}|R_dU x| \ge
\tau \right).
\end{equation}
This will follow from the Fubini theorem once we show
a similar inequality for {\em conditional }
probabilities for any fixed $U$.  To that end, consider any
$U \in O(m)$ verifying the condition on the right hand side of
(\ref{polar}),  i.e., such that there exists $y \in R_d U S$ for which
$|y| \ge \tau$. Then
\begin{eqnarray*}
\pp \left(\max_{x \in S}|(A A^*)^{1/2}R_dU x|  \ge
\tau \right)  & \ge &
\pp \left(|(A A^*)^{1/2}y|
\ge \tau \right)  \\
&  = &
\pp \left(|A^*y| \ge \tau \right)
\ge \pp \left(|A^*y| \ge |y| \right).
\end{eqnarray*}
Again by rotational invariance of the Gaussian measure, the
distribution of the random vector $A^*y$ depends only on
$|y|$ and is an appropriate multiple of the standard Gaussian
vector on $\R^m$.  Accordingly,
$\pp \left(|A^*y| \ge s  \right)$
may be expressed in terms of the appropriate $\chi^2$ or Gamma
distributions.  For example, for any $y \neq 0$,
$$
\pp \left(|A^*y| \ge |y| \right) = \pp \left(\chi^2_m \ge m
\right) = \left(\Gamma(m/2)\right)^{-1} \int_{m/2}^\infty
e^{-u} u^{m/2-1} du  =: \gamma_m .
$$
It follows from the central limit theorem that $\gamma_m \ra 1/2$ as
$m \ra \infty$. Moreover, it is easily verified that  $\gamma_m \ge
\gamma_2 = e^{-1}$ if $m \ge 2$.  This shows  (\ref{polar})
for $m > 1$ and concludes the proof of the Lemma.  \qed

\footnotesize
\address

\end{document}